J. Mikesh, S. E. Stepanov

# Projective Maps from the Perspective of Elliptic Differential Operators

**Abstract.** This paper develops an analytical approach to the study of the geometry of projective maps using the theory of elliptic differential operators. We construct two elliptic operators of second and fourth order, whose kernels characterize projective diffeomorphisms between Riemannian manifolds and one-parameter groups of projective diffeomorphisms (transformations) of a Riemannian manifold onto itself, respectively. This approach establishes a natural correspondence between analytical and geometric properties, enabling the study of projective diffeomorphisms via operator-theoretic methods. The proposed framework provides a new understanding of projective structures on Riemannian manifolds and extends classical results in differential geometry.

**Keywords:** Riemannian manifolds, elliptic operators, projective diffeomorphisms, spectral theory.

**MSC2010:** 53C20; 35J15; 35J30

## §1. Introduction

Differential operators are often employed as a primary tool for investigating geometric and spectral properties of manifolds (see [1]–[6]). Examples include various types of Laplacians (cf. [7, pp. 77–83]; [9]; [10]), the Dirac operator and its analogues (cf. [2, pp. 92–94]; [7, p. 81]; [8, pp. 363–722]; [11]). Elliptic operators (see, e.g., [2]) have recently become an integral part of modern Riemannian geometry, providing powerful methods for analyzing and understanding the structure of manifolds (see, e.g., [3]; [7]; [12]; [13]).

Within the context of Riemannian geometry, projective mappings are diffeomorphisms that preserve geodesics (see, for example, [14]; [15]; [16]). The central concept here is that of projective equivalence of metrics: two Riemannian metrics on a manifold are said to be projectively equivalent if they share the same geodesics. A series of works (see [17]; [18]; [19]) have applied differential operators, particularly the Nijenhuis operator, to the study of projectively equivalent metrics. However, these studies did not fully situate projective diffeomorphisms within the general framework of differential operator theory.

In contrast, the present work demonstrates how projective diffeomorphisms of Riemannian manifolds can be studied using the theory of elliptic operators. In particular, we show how the analytical properties of these operators are connected with one-parameter groups of projective diffeomorphisms acting on a Riemannian manifold. Specifically, we define and investigate two elliptic operators of second and fourth order, which respectively characterize projective diffeomorphisms between Riemannian manifolds and one-parameter groups of local diffeomorphisms of a Riemannian manifold onto itself.

The proposed approach — studying the geometry of projective diffeomorphisms and their one-parameter groups through elliptic differential operator theory — is novel and provides a fresh analytical perspective on classical problems in Riemannian geometry.

## §2. Preliminaries

Let $(M, g)$ be a connected Riemannian manifold of dimension $n \geq 2$ with metric tensor $g$ and Levi-Civita connection $\nabla$. Denote by $TM$ and $T^*M$ its tangent and cotangent bundles, respectively, and by $\otimes^r T^*M$, $S^r M := S^r(T^*M)$, and $\Lambda^r M := \Lambda^r(T^*M)$ the bundles of covariant $r$-tensors and their symmetric and antisymmetric subbundles. We write $C^\infty(TM)$, $C^\infty(T^*M)$, $C^\infty(\otimes^r T^*M)$, $C^\infty(S^r M)$, and

$C^\infty(\Lambda^r M)$ for the vector spaces of smooth sections of these bundles. In particular, $C^\infty(T^*M)$ will be denoted $\Omega^1(M)$ throughout this paper.

The metric on each of these bundles is induced by the Riemannian metric $g$, which we denote by the same symbol. On a compact boundaryless manifold $M$, the Hodge metric defines an $L^2$-inner product on the space of sections $C^\infty(\otimes^r T^*M)$ by

$$\langle \omega, \theta \rangle_{L^2} = \int_M g(\omega, \theta)\, dv_g,$$

for all $\omega, \theta \in C^\infty(\otimes^r T^*M)$. This equips the space of square-integrable functions or tensor fields on $M$ with the structure of a Hilbert space with respect to this inner product, providing a natural setting for the study of differential operators such as Laplace-type operators and their generalizations.

Given a linear differential operator $D: C^\infty(F) \to C^\infty(E)$ acting between spaces of smooth sections of subbundles $F$ and $E$ of $\otimes^r T^*M$ on a compact boundaryless Riemannian manifold $(M, g)$, the Hodge inner product allows one to define the formal adjoint $D^*$ via (see [2, p. 70]; [7, p. 626])

$$\langle D\omega, \eta \rangle_{L^2} = \langle \omega, D^*\eta \rangle_{L^2}, \quad \forall \omega \in C^\infty(F), \eta\, C^\infty(E).$$

As an example, the operator $\nabla^*$, formally adjoint to the covariant derivative $\nabla: C^\infty(\otimes^r T^*M) \to C^\infty(T^*M \otimes \otimes^r T^*M)$, is given by

$$(\nabla^*\omega')(Y_1, \ldots, Y_r) = -\sum_{i=1}^n (\nabla_{e_i}\omega')(e_i, Y_1, \ldots, Y_r),$$

where $\omega' \in C^\infty(T^*M \otimes \otimes^r T^*M)$, $Y_1, \ldots, Y_r \in C^\infty(TM)$, and $\{e_1, \ldots, e_n\}$ is an arbitrary orthonormal basis of $T_xM$ for each $x \in M$ (cf. [7, pp. 34–35]).

As a second example, consider the symmetric differentiation operator $\delta^*: C^\infty(S^r M) \to C^\infty(S^{r+1}M)$ (see [7, p. 514]). On the space of 1-forms $\Omega^1(M)$, this operator is expressed as

$$\delta^*\theta = \frac{1}{2}\mathcal{L}_\xi g,$$

where $\mathcal{L}_\xi g$ denotes the Lie derivative of the metric tensor $g$ along the vector field $\xi := \theta^\#$, the metric dual of $\theta$ (see [7, p. 55]; [20]). In particular, $\delta^*\theta = 0$ if and only if $\xi = \theta^\#$ is a Killing vector field, i.e., an infinitesimal isometry (cf. [7, p. 55]). Every infinitesimal isometry $\xi$ generates, locally near any point $x \in M$, a one-parameter group of local diffeomorphisms preserving distances on $(M, g)$.

The formal adjoint of $\delta^*$ with respect to the Hodge metric is the divergence operator

$$\delta: C^\infty(S^{r+1}M) \to C^\infty(S^r M),$$

which is the restriction of $\nabla^*$ to the symmetric subbundle $S^{r+1}M$ (see [7, p. 514]). If $\delta\theta = 0$ for $\theta \in \Omega^1(M)$, the form is called co-closed and the vector field $\xi = \theta^\#$ is solenoidal. In this context, the classical Berger–Ebin $L^2$-orthogonal decomposition holds (cf. [20]):

$$C^\infty(S^2 M) = \ker \delta \oplus_{L^2} \operatorname{im} \delta^*.$$

Since the operators considered later in this paper are of the Stein–Weiss type $D^*D$, we describe this construction in detail. Let $D: C^\infty(F) \to C^\infty(E)$ be a linear differential operator as above, with formal adjoint $D^*$. If $D$ is injective, one has the $L^2$-orthogonal decomposition (see [20])

$$C^\infty(E) = \ker D^* \oplus_{L^2} \operatorname{im} D.$$

The operator $D^*D$ is elliptic, its kernel $\ker(D^*D)$ is finite-dimensional, and the $L^2$-orthogonal decomposition

$$C^\infty(F) = \ker(D^*D) \oplus_{L^2} \operatorname{im}(D^*D)$$

holds (cf. [7, p. 632]). Examples of elliptic operators of the form $D^*D$ include the rough Laplacian (Bochner Laplacian) $\nabla^*\nabla: C^\infty(\otimes^r T^*M) \to C^\infty(\otimes^r T^*M)$ (cf. [7, p. 77]; [21]), the Tachibana operator $T^*T: C^\infty(\Lambda^r M) \to C^\infty(T^*M \otimes \Lambda^r M)$ (cf. [22,

23]), and its special case, the Alfors Laplacian $S^*S: C^\infty(T^*M) \to C^\infty(S_0^2 M)$ acting on the bundle of trace-free symmetric 2-tensors $S_0^2 M$(cf. [21, 24]).

Well-known examples of elliptic operators of the type $D^*D$ are the Stein–Weiss operators, which have been studied in the papers [11, 21] and in the monograph [22, pp. 260–279]. It has been shown that ellipticity for operators of the form $D^*D$ can be achieved for generalized gradients $D$, i.e., first-order equivariant differential operators acting between irreducible vector bundles with structure groups $O(n)$, $SO(n)$, or $\mathrm{Spin}(n)$. The monograph [22, pp. 260–279] provides numerous examples of Stein–Weiss operators $D^*D$, and in [21] the spectral decomposition of Stein–Weiss elliptic operators $D^*D$ on the standard sphere $\mathbb{S}^n$ was obtained for each generalized gradient $D$.

## §3. Definition and Study of the Sinyukov Elliptic Operator in the Theory of Projective Mappings

Let $(M, g)$ be a connected compact Riemannian manifold of dimension $n \geq 2$. Consider the first-order linear differential operator

$$S: C^\infty(M, S^2(T^*M)) \to C^\infty(M, T^*M \otimes S^2(T^*M))$$

defined for a symmetric 2-tensor field $\phi$ by

$$(S\phi)(X, Y, Z) = (\nabla_X \phi)(Y, Z) - \frac{1}{n+1}(g(X, Y)\,\mathrm{div}\,\phi(Z) + g(X, Z)\,\mathrm{div}\,\phi(Y)),$$

where $\nabla$ is the Levi-Civita connection of $g$ and $\mathrm{div}\,\phi(W) := \mathrm{trace}_g(\nabla \phi)(\cdot, W)$.

The kernel $\ker S$ is characterized by

$$(\nabla_X \phi)(Y, Z) = \frac{1}{n+1}(g(X, Y)\mathrm{div}\,\phi(Z) + g(X, Z)\mathrm{div}\,\phi(Y)).$$

A solution is called *trivial* if $\phi = Cg$ for some constant $C$.

A diffeomorphism $f: (M, g) \to (\bar{M}, \bar{g})$ is called *projective* if it maps every geodesic of $(M, g)$ to a geodesic of $(\bar{M}, \bar{g})$ (cf. [15,26,27]). Sinyukov [28, Thm. 1] showed that $(M, g)$ admits a nontrivial projective mapping if and only if there exists $\phi \in C^\infty(M, S^2(T^*M))$ satisfying

$$\nabla_X \varphi = g(X, \cdot) \otimes \omega + g(\cdot, \cdot) \otimes \omega(X)$$

for some 1-form $\omega \in \Omega^1(M)$. Locally, the projectively related metric $\bar{g}$ is given by

$$\bar{g}_{ij} = e^{2\rho} \varphi^{kl} g_{ki} g_{lj}, \quad \nabla_i \rho = -\omega_k \varphi^{kl} g_{li}, \quad \nabla_i(g^{kl} \varphi_{kl}) = 2\omega_i,$$

and it follows that $\omega = \frac{1}{n+1} \operatorname{div} \varphi$. Thus, $\ker S$ corresponds to pseudo-Riemannian metrics projectively diffeomorphic to $(M, g)$.

The family of manifolds $(\bar{M}, \bar{g})$ admitting projective maps from $(M, g)$ depends on at most $\frac{1}{2}(n+1)(n+2)$ parameters, with equality iff $(M, g)$ has constant curvature (cf. [28, Thm. 2]).

**Lemma 1.** *Let $(M, g)$ be compact and $S$ as above. Then the operator*

$$S^* : C^\infty(M, T^*M \otimes S_0^2(TM)) \to C^\infty(M, S^2(T^*M))$$

*defined by*

$$(S^*\Phi)(Y, Z) = -\operatorname{trace}_g(\nabla \Phi)(\cdot, Y, Z) + \frac{1}{n+1}((\nabla_Y \operatorname{tr}_g \Phi)(Z) + (\nabla_Z \operatorname{tr}_g \Phi)(Y))$$

*is formally adjoint to $S$ with respect to the $L^2$-inner product.*

**Proof.** We assume that the manifold $M$ is connected and compact (without boundary). All subsequent arguments are carried out in local coordinates $x^1, \ldots, x^n$ associated with an arbitrary chart $(U, \psi)$ of the manifold $M$.

The operator formally adjoint to $S$ is the operator

$$S^* : C^\infty(T^*M \otimes S_0^2 M) \to C^\infty(S^2 M),$$

defined by the identity

$$\langle S\varphi, \Phi \rangle_{L^2} = \langle \varphi, S^* \Phi \rangle_{L^2}$$

for arbitrary $\varphi \in C^\infty(S^2 M)$ and $\Phi \in C^\infty(T^*M \otimes S_0^2 M)$.

By integration by parts, we obtain

$$\langle S\varphi, \Phi \rangle_{L^2} = \int_M (S\varphi)_{kij} \Phi^{kij} dv_g$$

$$= \int_M \left( \nabla_k \varphi_{ij} - \frac{1}{n+1}(g_{ki}(\nabla^l \varphi_{lj}) + g_{kj}(\nabla^l \varphi_{li})) \right) \Phi^{kij} dv_g$$

$$= \int_M (\nabla_k \varphi_{ij}) \Phi^{kij} dv_g - \frac{1}{n+1} \int_M (g_{ki}(\nabla^l \varphi_{lj})) \Phi^{kij} dv_g$$

$$-\frac{1}{n+1}\int_M (g_{kj}(\nabla^l \varphi_{li}))\, \Phi^{kij}\, dv_g$$

$$=\int_M \varphi_{ij}\left(-\nabla_k \Phi^{kij} + \frac{1}{n+1}\left(\nabla^i \Phi_k{}^{kj} + \nabla^j \Phi_k{}^{ki}\right)\right) dv_g = \langle \varphi, S^*\Phi\rangle_{L^2}.$$

where $\Phi_{kij}$ denote the local components of the tensor field $\Phi \in C^\infty(T^*M \otimes S_0^2 M)$.

Thus, the equality to be proved follows.

Based on the lemma, we prove that the following theorem is valid.

**Theorem 1.** *There is an $L^2$-orthogonal decomposition*

$$C^\infty(M, T^*M \otimes S^2(T^*M)) = \ker S \oplus_{L^2} \operatorname{im} S^*,$$

*where $S$ is injective and $S^*$ is the formal adjoint.*

**Proof.** The principal symbol of operator (3.1) is given by

$$(\sigma_S(\theta)\varphi)_{kij} = \theta_k \varphi_{ij} - \frac{1}{n+1}(g_{ki}\,\theta^l \varphi_{lj} + g_{kj}\,\theta^l \varphi_{li}),$$

for $\theta \in T_x^*M$.

This follows from the well-known fact that for the covariant derivative $\nabla$ and its formal adjoint $\nabla^*$, the principal symbols are respectively (see, e.g., [7, p. 628]; [2, pp. 79, 87])

$$(\sigma_\nabla(\theta)\varphi)_{kij} = i\,\theta_k \varphi_{ij}, \quad (\sigma_{\nabla^*}(\theta)\varphi)_j = -i\,\theta^k \varphi_{kj}.$$

In the formula for the principal symbol we omit, as is customary (see [2, pp. 79, 87]), the common factor $i$ and retain only the algebraic part, since this factor does not affect injectivity.

Suppose now that $\sigma_S(\theta)\varphi = 0$, i.e.,

$$\theta_k \varphi_{ij} = \frac{1}{n+1}(g_{ki}\,\theta^l \varphi_{lj} + g_{kj}\,\theta^l \varphi_{li}).$$

Denote by $\alpha_j$ the contraction $\theta^l \varphi_{lj}$. Then the above identity can be rewritten as

$$\theta_k \varphi_{ij} = \frac{1}{n+1}(g_{ki}\,\alpha_j + g_{kj}\,\alpha_i).$$

Contracting both sides of the equality obtained above with $\theta^j \theta^k \alpha^i$, we obtain

$$(\varphi_{lj}\theta^l \theta^j)^2 = n\,|\theta|^2\,|\alpha|^2.$$

On the other hand, applying the Cauchy–Schwarz inequality to the vectors $\alpha$ and $\theta$, we have

$$|\alpha \cdot \theta| \leq |\theta| \ |\alpha|,$$

where

$$\alpha \cdot \theta := \alpha_i \theta^i = \varphi_{lj} \theta^l \theta^j.$$

Consequently,

$$(\varphi_{lj}\theta^l\theta^j)^2 \leq |\theta|^2 \ |\alpha|^2.$$

Combining the two relations yields

$$n \ |\theta|^2 \ |\alpha|^2 \leq |\theta|^2 \ |\alpha|^2.$$

If $|\theta|^2 \neq 0$ and $|\alpha|^2 \neq 0$, this implies $n \leq 1$, which contradicts the assumption $n \geq 2$. Hence, for $n \geq 2$ we must have $|\alpha|^2 = 0$, that is, $\alpha = 0$.

Substituting $\alpha = 0$ into

$$\theta_k \varphi_{ij} = \frac{1}{n+1}(g_{ki}\,\alpha_j + g_{kj}\,\alpha_i),$$

we obtain

$$\theta_k \varphi_{ij} = 0 \ \text{ for all } k, i, j.$$

If $\theta \neq 0$, then there exists an index $k$ such that $\theta_k \neq 0$, and therefore the above relations imply $\varphi_{ij} = 0$.

Thus, for every nonzero $\theta$ and $n \geq 2$, the condition $\sigma_S(\theta)\varphi = 0$ implies $\varphi = 0$. Hence, for $n \geq 2$, the principal symbol $\sigma_S(\theta)$ is injective for all $\theta \neq 0$, as required.

In this case, the orthogonal decomposition

$$C^\infty(T^*M \otimes S^2 M) = \ker S \oplus_{L^2} \text{im } S^*$$

holds with respect to the Hodge metric (see [20, Corollary 4.2]).

Let us now prove the corollary of this theorem.

**Corollary 1.** *The second-order operator*

$$S^*S: C^\infty(M, S^2(T^*M)) \to C^\infty(M, S^2(T^*M))$$

*defined by*

$$(S^*S\,\varphi)(Y, Z) = \bar{\Delta}\varphi(Y, Z) + \frac{1}{n+1}((\nabla_Y \text{div } \varphi)(Z) + (\nabla_Z \text{div } \varphi)(Y))$$

*is elliptic, self-adjoint, nonnegative, and has finite-dimensional kernel* $\ker S^*S = \ker S \subset C^\infty(M, S^2(T^*M))$.

**Proof.** Direct computations yield the above identities, which define the second-order differential operator

$$S^*S(\varphi) := S^*(S\varphi),$$

acting on the space of symmetric 2-tensors $C^\infty(S^2 M)$. Since, as established in Theorem 1, the operator $S$ has injective principal symbol, it follows that the operator $S^*S$ is elliptic (see [7, p. 629]; [20, p. 383]). Consequently, its kernel ker $S^*S$ is a finite-dimensional vector space contained in $C^\infty$ (see [7, p. 631]). The self-adjointness and non-negativity of the elliptic operator $S^*S$ follow from the identities (see, for example, [7, pp. 631–632])

$$\langle S^*S\,\omega, \omega'\rangle_{L^2} = \langle S\omega, S\omega'\rangle_{L^2} = \langle \omega, S^*S\,\omega'\rangle_{L^2},$$

where $\langle S\varphi, S\varphi\rangle_{L^2} = \|S\varphi\|^2_{L^2} \geq 0$. This completes the proof.

On the basis of the standard properties of elliptic operators on compact (boundaryless) Riemannian manifolds (see, for example, [4], [5], [6], and [7, p. 632]), we conclude that the kernel ker $(S^*S)$ coincides with the kernel of the operator $S$. Therefore, ker $(S^*S)$ consists precisely of all (pseudo-)Riemannian metrics $\bar{g}$ that are projectively diffeomorphic to $(M, g)$, just as does the kernel of $S$. Moreover, ker$(S^*S)$ is finite-dimensional and is contained in $C^\infty$ (see [7, p. 631]). Consequently, the kernel of $S$ is also finite-dimensional, which is consistent with Sinjukov's result cited above concerning the dimension of the space of solutions to equations (3.3).

All of the above justifies calling $S^*S$ the *Sinjukov operator*. The following corollary holds.

**Corollary 2.** *The Sinyukov operator $S^*S$ has the kernel consisting of pseudo-Riemannian metrics projectively diffeomorphic to $(M,g)$, and there is the $L^2$-orthogonal decomposition*

$$C^\infty(M, S^2(T^*M)) = \operatorname{im}(S^*S) \oplus_{L^2} \ker(S^*S).$$

**Remark.** In the monograph [22, pp. 265–266], a decomposition of the space

$$T_x^*M \otimes S^2(T_x^*M)$$

into an orthogonal direct sum of three subspaces was obtained, two of which are pointwise $O(n)$-irreducible. As a consequence, there exist only two generalized gradients

$$D: C^\infty(S^2 M) \to C^\infty(T^*M \otimes S^2 M),$$

and, based on them, one can construct two Stein–Weiss elliptic operators of the form $D^*D$.

The differential operator

$$S: C^\infty(S^2 M) \to C^\infty(T^*M \otimes S^2 M)$$

is not among the generalized gradients identified there; therefore, the Sinjukov operator $S^*S$ is not a Stein–Weiss elliptic operator.

The following statement describes the structure of the kernel of the Sinjukov operator $S^*S$ on a compact Riemannian manifold $(M, g)$ of nonpositive sectional curvature.

**Corollary 3.** *If $(M, g)$ has non-positive sectional curvature, $\ker(S^*S)$ consists of parallel 2-tensors. If the curvature is strictly negative in all 2-planes at some point, then $\ker(S^*S) = \{Cg \mid C \in \mathbb{R}\}$.*

**Proof.** By Corollary 2, on a compact manifold $(M, g)$ the kernel of the first-order linear differential operator

$$S: C^\infty(S^2 M) \to C^\infty(T^*M \otimes S_0^2 M),$$

acting according to formula (3.1), coincides with the kernel of the elliptic operator $S^*S$ defined by (3.6). Hence, if $\varphi \in C^\infty(S^2 M)$ lies in $\ker(S^*S)$, then it satisfies

$$\nabla_k \varphi_{ij} = \frac{1}{n+1}(g_{ki} \nabla^l \varphi_{lj} + g_{kj} \nabla^l \varphi_{li}).$$

On the other hand, if $(M, g)$ is a connected compact (boundaryless) Riemannian manifold, then for the local components $\varphi_{ij}$ of an arbitrary tensor $\varphi \in C^\infty(S^2 M)$ the following integral identity holds (see [30]):

$$\int_M (K(\varphi, \varphi) + \nabla^k \varphi^{ij} \nabla_i \varphi_{kj} - \nabla_i \varphi^{ij} \nabla^k \varphi_{kj}) \, dv_g = 0. \qquad (3.6)$$

Here, by [6, p. 592] and [20],

$$K(\varphi,\varphi) := \sum_{i<j} \sec(e_i, e_j)(\rho_i - \rho_j)^2,$$

where $\{e_1, \ldots, e_n\}$ is an orthonormal basis of $T_x M$ at an arbitrary point $x \in M$ chosen so that $\varphi(e_i, e_j) = \rho_i \delta_{ij}$, and $\sec(e_i, e_j)$ denotes the sectional curvature of the 2-plane $\sigma = \text{span}\{e_i, e_j\} \subset T_x M$.

If $\sec \leq 0$, then $K(\varphi, \varphi) \leq 0$; moreover, $K(\varphi, \varphi) < 0$ whenever $\sec(e_i, e_j) < 0$ for every pair $i \neq j$.

For an arbitrary symmetric 2-tensor $\varphi \in C^\infty(S^2 M) \cap \ker(S^*S)$, identity (3.6) can be rewritten as

$$\int_M \left( \sum_{i<j} \sec(e_i, e_j)(\rho_i - \rho_j)^2 - \frac{(n-1)(n+2)}{4(n+1)^2} \parallel \delta\varphi \parallel^2 \right) dv_g = 0,$$

where $\parallel \delta\varphi \parallel^2 := g(\delta\varphi, \delta\varphi) = \nabla_i \varphi^{ij} \nabla^k \varphi_{kj}$.

It is now clear that, when $\sec \leq 0$ and $\varphi \in C^\infty(S^2 M) \cap \ker(S^*S)$, such a tensor must be covariantly constant, i.e., a parallel symmetric 2-tensor.

If, moreover, the sectional curvature at some point $x \in M$ is negative in all 2-planes in $T_x M$, then the only tensor $\varphi \in C^\infty(S^2 M) \cap \ker(S^*S)$ on the compact manifold $(M, g)$ is $\varphi = C g$, $C = \text{const}$.

**Remark.** The above-mentioned result of Sinjukov [27] on the finiteness of the dimension of the space of solutions to Sinjukov equations is local. This is consistent with Corollary 3 proved here, since it excludes the compact hyperbolic space $\mathbb{H}^n$ from consideration.

## §4. Spectral properties of the Sinyukov operator

Let us briefly discuss the basic spectral properties of the **Sinjukov operator** $S^*S$. According to the theory of elliptic operators on a compact Riemannian manifold $(M, g)$, a self-adjoint nonnegative elliptic operator such as $S^*S$ has discrete spectrum

$$\text{Spec}(S^*S) = \{0 < \mu_1 \leq \mu_2 \leq \cdots \to \infty\}.$$

The eigenvalue 0 corresponds to tensors in ker $(S^*S)$, i.e., to the space of solutions of the equation $S^*S(\varphi) = 0$. As established above, ker$(S^*S)$ is finite-dimensional; therefore, the multiplicity of the zero eigenvalue equals

$$\dim_{\mathbb{R}} \ker (S^*S) \leq \frac{1}{2}(n+1)(n+2).$$

The positivity of the principal symbol of $S^*S$ implies the estimate (see [29]): there exist constants $c > 0$ and $C \geq 0$ such that for all $\varphi \in C^\infty(S^2M)$,

$$\langle S^*S\varphi, \varphi \rangle_{L^2} = \| S\varphi \|^2_{L^2} \geq c \| \nabla\varphi \|^2_{L^2} - C \| \varphi \|^2_{L^2},$$

where $c > 0$ depends only on the leading coefficients of the operator. On the $L^2$-orthogonal complement of ker $S$, the spectrum of $S^*S$ admits a strictly positive lower bound of the form

$$\langle S^*S\varphi, \varphi \rangle_{L^2} \geq \alpha \| \varphi \|^2_{L^2}$$

for some constant $\alpha > 0$.

We now turn to an explicit estimate of the eigenvalues of $S^*S$ in the model case of the standard sphere $(S^n, g_{\text{can}})$, where $g_{\text{can}}$ denotes the canonical metric of constant sectional curvature 1.

Recall the $L^2$-orthogonal decomposition (see [7, p. 130]):

$$C^\infty(S^2M) = (\text{Im } \delta^* + C^\infty(M) \cdot g) \oplus_{L^2} (\delta^{-1}(0) \cap \text{trace}_g^{-1}(0)),$$

where $\delta^{-1}(0) \cap \text{trace}_g^{-1}(0)$ is the space of trace-free, divergence-free symmetric 2-tensors, the so-called *TT-tensors*.

To estimate the eigenvalues of $S^*S$ on $(\mathbb{S}^n, g_{\text{can}})$, we use the known spectra of the scalar, vector, and tensor Laplacians, together with the relations between $S^*S$, $\nabla^*\nabla$, and the Lichnerowicz Laplacian $\Delta_L$ (see [7, pp. 80, 178]) on each component of the decomposition (4.1).

1. *The trace component $C^\infty(M) \cdot g$*

In this case, the eigenvalues $\mu_k$ of $S^*S$ are expressed in terms of the eigenvalues of the scalar Laplacian $\bar{\Delta}$, assuming sectional curvature equal to 1. As a result, the first nonzero eigenvalue equals $\mu = n$.

2. *The component Im $\delta^*$*

Here, the eigenvalues $\mu_k$ of $S^*S$ are expressed in terms of the spectrum of the vector Laplacian under the condition $\text{Ric} = (n-1)g$. For the part corresponding to gradients, the spectrum depends directly on the eigenvalues of the vector Laplacian:

$$\mu_k = \frac{k(n+k-1)}{n}, k = 1,2,\ldots,$$

where $k$ denotes the order of the vector spherical harmonics.

3. *The TT-component $\delta^{-1}(0) \cap trace_g^{-1}(0)$.*

In this case, the eigenvalues $\mu_k$ of $S^*S$ coincide with the spectrum of the Lichnerowicz Laplacian $\Delta_L$ on tensor spherical harmonics (see details in [31]):

$$\mu_k = \frac{k(n+k-1)}{n}, k = 1,2,\ldots,$$

where $k$ corresponds to the degree of the tensor harmonic.

Thus, on the standard sphere, the spectrum of the Sinjukov operator $S^*S$ can be computed explicitly on each irreducible component of the above $L^2$-orthogonal decomposition, and it is determined by the well-known spectra of the scalar, vector, and tensor Laplacians.

## §5. Definition and Study of the Eisenhart Elliptic Operator in the Theory of Projective Transformations

Consider the second-order differential operator

$$E: \Omega^1(M) \to C^\infty(T^*M \otimes S_0^2 M),$$

which in a local coordinate system $x^1,\ldots,x^n$ of an arbitrary chart of the manifold $M$ is given by

$$(E\theta)_{kij} = 2(n+1)\nabla_k(\delta^*\theta)_{ij} + 2g_{ij}\nabla_k(\delta\theta) + g_{kj}\nabla_i(\delta\theta) + g_{ki}\nabla_j(\delta\theta),$$

for any 1-form $\theta \in \Omega^1(M)$.

The kernel $\ker E$ of the operator $E$ is determined by the *Eisenhart equations* (see [14, p. 275]; [32])

$$2(n+1)\nabla_k(\delta^*\theta)_{ij} = -2g_{ij}\nabla_k(\delta\theta) - g_{kj}\nabla_i(\delta\theta) - g_{ki}\nabla_j(\delta\theta).$$

Each solution of these equations corresponds to a smooth vector field $\xi = \theta^{\sharp}$, which generates in a neighborhood of any point $x \in M$ a local one-parameter group of diffeomorphisms of the Riemannian manifold $(M, g)$ that map geodesics to geodesics (see [14, pp. 274–277]). Such a vector field $\xi$ is called an infinitesimal projective transformation of $(M, g)$ (see [33]), with a trivial example being an infinitesimal isometry.

Analysis of the integrability conditions of $E$ shows (see [33]) that the set of all infinitesimal projective diffeomorphisms of a connected manifold $(M, g)$ forms a finite-dimensional Lie algebra whose dimension does not exceed $n(n + 2)$, with equality if and only if $(M, g)$ has constant curvature.

For a connected compact Riemannian manifold without boundary $(M, g)$, one can define the formal adjoint operator $E^*$ with respect to the Hodge metric:

**Lemma 2.** *Let $(M, g)$ be a connected compact Riemannian manifold without boundary of dimension $n \geq 2$, and let $E: \Omega^1(M) \to C^\infty(T^*M \otimes S_0^2 M)$ be the second-order differential operator above. Then the operator*

$$(E^*\Phi)_k = 2(n+1)\nabla^i \nabla^j \Phi_{jik} + 2\nabla_k \nabla_i \Phi_j^{ji},$$

*for any tensor field $\Phi \in C^\infty(T^*M \otimes S_0^2 M)$ with local components $\Phi_{kij}$, is formally adjoint to $E$ with respect to the Hodge metric.*

**Proof.** Let $(M, g)$ be a connected compact (boundaryless) Riemannian manifold. We work in local coordinates $x^1, \ldots, x^n$ associated with an arbitrary chart $(U, \psi)$ on $M$.

In the present situation, by the definition of the formal adjoint we have

$$\langle E\theta, \Phi \rangle_{L^2} = \langle \theta, E^*\Phi \rangle_{L^2}$$

for arbitrary $\Phi \in C^\infty(T^*M \otimes S_0^2 M)$ and $\theta \in \Omega^1(M)$. Hence,

$$\langle E\theta, \Phi\rangle_{L^2} = \int_M (2(n+1)\,\nabla_k(\delta^*\theta)_{ij} + g_{kj}\nabla_i(\delta\theta) + g_{ki}\nabla_j(\delta\theta))\,\Phi^{kij}\,dv_g$$

$$= \int_M \theta_k\,(E^*\Phi)^k\,dv_g = \langle \theta, E^*\Phi\rangle_{L^2}.$$

In the left-hand side, we integrate the first term by parts and obtain

$$2(n+1)\int_M \nabla_k\,(\delta^*\theta)_{ij}\,\Phi^{kij}\,dv_g = -2(n+1)\int_M (\delta^*\theta)_{ij}\,\nabla_k\Phi^{kij}\,dv_g$$

$$= 2(n+1)\int_M \theta_k\,(\nabla_i\nabla_j\Phi^{jik})\,dv_g.$$

Next, we integrate by parts the second term:

$$\int_M (g_{kj}\nabla_i(\delta\theta) + g_{ki}\nabla_j(\delta\theta))\Phi^{kij}\,dv_g$$

$$= \int_M (\nabla_i(\delta\theta)\,\Phi_j{}^{ij} + \nabla_j(\delta\theta)\,\Phi_i{}^{ij})dv_g = 2\int_M (\nabla_i(\delta\theta)\,\Phi_j{}^{ij})dv_g =$$

$$= -2\int_M (\delta\theta)\,\nabla_i\Phi_k{}^{ki}\,dv_g = 2\int_M \theta^k\,(\nabla_k\nabla_i\Phi_j{}^{ji})dv_g$$

where $\Phi_k{}^{kj} = g^{ik}g^{jm}\Phi_{ikm}$.

Collecting the contributions, we obtain

$$(E^*\Phi)_k = 2(n+1)\,\nabla^i\nabla^j\Phi_{jik} + 2\nabla_k\nabla_i\Phi_j{}^{ji}$$

where $\nabla^l = g^{lk}\nabla_k$.

**Theorem 2.** *Let $(M,g)$ be a connected compact Riemannian manifold without boundary of dimension $n \geq 2$. Then there exists an $L^2$-orthogonal decomposition*

$$C^\infty(T^*M) = \ker E \oplus_{L^2} \operatorname{im} E^*,$$

*where $E$ is the injective second-order operator above, and $E^*$ is its formal adjoint.*

**Proof.** We rewrite identity

$$(E\theta)_{kij} = 2(n+1)\nabla_k(\delta^*\theta)_{ij} + 2g_{ij}\nabla_k(\delta\theta) + g_{kj}\nabla_i(\delta\theta) + g_{ki}\nabla_j(\delta\theta),$$

in the following form:

$$(E\theta)_{kij} = (n+1)\,\nabla_k(\nabla_i\theta_j + \nabla_j\theta_i) - 2\,g_{ij}\,\nabla_k(\nabla^l\theta_l) - g_{kj}\,\nabla_i(\nabla^l\theta_l)$$
$$- g_{ki}\,\nabla_j(\nabla^l\theta_l),$$

since $(\delta^*\theta)_{ij} := \frac{1}{2}(\nabla_i\theta_j + \nabla_j\theta_i)$ and $\delta\theta := -\nabla^l\theta_l$ for the local components $\theta_i$ of the 1-form $\theta$.

The principal symbols of the operators $\delta^*$ and $\delta$ are known:

$$\sigma_{\delta^*}(\omega)(\theta)_{ij} = \frac{1}{2}(\omega_i\theta_j + \omega_j\theta_i), \quad \sigma_\delta(\omega)(\theta) = -\omega^j\theta_j.$$

Therefore, for a nonzero covector $\omega \in T_x^*M$, the principal symbol $\sigma_E(\omega)$ of the operator $E$ acts by

$$\sigma_E(\omega): \theta \mapsto i((n+1)\,\omega_k(\omega_i\theta_j + \omega_j\theta_i) - 2\,g_{ij}\,\omega_k(\omega^l\theta_l) - g_{kj}\,\omega_i(\omega^l\theta_l)$$
$$- g_{ki}\,\omega_j(\omega^l\theta_l)).$$

One can show that, for every point $x \in M$ and every covector $\omega \neq 0$ in $T_x^*M$, the linear map

$$\sigma_E(\omega): T_x^*M \to T_x^*M \otimes S_0^2(T_x^*M)$$

has trivial kernel: if $\sigma_E(\omega)\theta = 0$, then $\theta = 0$.

Indeed, contracting the last identity in the indices $k$ and $i$, and using $\omega \neq 0$, one successively derives first that $\omega^l\theta_l = 0$, and then that $\omega_i\theta_j + \omega_j\theta_i = 0$. For $\omega \neq 0$, these two conditions imply $\theta = 0$.

Thus, the principal symbol $\sigma_E(\omega)$ is injective for every $\omega \neq 0$. Consequently, on a compact (boundaryless) manifold $(M, g)$ the operator $E$ is overdetermined elliptic (see [7, p. 629]) or, equivalently, injectively elliptic (see [20]). In this case, the orthogonal decomposition

$$C^\infty(T^*M) = \ker E \oplus_{L^2} \operatorname{im} E^*$$

holds with respect to the Hodge metric (see [13, Corollary 4.2]).

**Remark.** In the monograph [22, p. 277], in the course of determining generalized gradients $D$ for the subsequent construction of elliptic operators of the form

$$D^*D: C^\infty(S_0^2 M) \to C^\infty(T^*M \otimes S_0^2(T^*M))$$

and of "strong Laplacians" acting on $C^\infty(S_0^2 M)$, a decomposition of the space

$$T_x^* M \otimes S_0^2(T_x^* M)$$

(at an arbitrary point $x \in M$) into an orthogonal direct sum of three pointwise $O(n)$-irreducible components was proved. In Theorems 1 and 2 of the present work, this issue is resolved in the context of "geometry as a whole": namely, we obtain decompositions of the spaces $C^\infty(T^*M \otimes S^2 M)$ and $C^\infty(T^*M)$ into two $L^2$-orthogonal components in each of the two cases, respectively.

The principal symbol of $E$ is injective for every nonzero covector, implying that $E$ is overdetermined elliptic on compact $(M, g)$ without boundary (see [7, p. 629]; [20]). Hence, the $L^2$-orthogonal decomposition holds (see [13, Corollary 4.2]).

Consider the fourth-order operator

$$E^*E: \Omega^1(M) \to \Omega^1(M),$$

defined on $\Omega^1(M)$ by

$$E^*E(\theta) = E^*(E\theta)$$

for an arbitrary 1-form $\theta \in \Omega^1(M)$. Its properties are stated in the following corollary.

**Corollary 4.** *Under the same assumptions, the fourth-order operator*

$$E^*E: \Omega^1(M) \to \Omega^1(M), \quad E^*E(\theta) := E^*(E\theta),$$

*is elliptic, non-negative, self-adjoint, and has finite-dimensional kernel such that $\ker(E^*E) \subset C^\infty$.*

**Proof.** Since, as shown in the lemma, the operator $E$ has injective principal symbol, the operator $E^*E$ is elliptic (see [7, p. 629]; [20, p. 383]). Moreover, its kernel $\ker(E^*E)$ is a finite-dimensional vector space contained in $C^\infty$ (see [7, p. 631]).

The self-adjointness and non-negativity of the elliptic operator $E^*E$ follow from the identities (see, e.g., [7, pp. 631–632])

$$\langle E^*E\,\omega, \omega'\rangle_{L^2} = \langle E\omega, E\omega'\rangle_{L^2} = \langle \omega, E^*E\,\omega'\rangle_{L^2},$$

where

$$\langle E\varphi, E\varphi\rangle_{L^2} = \|E\varphi\|_{L^2}^2 \geq 0.$$

Note that in the absence of projective transformations on $(M, g)$, i.e. when $\ker E = \{0\}$, the operator $E^*E$ is strictly positive.

**Remark.** We present an alternative proof of the ellipticity of the operator $S^*S$. Namely, one can show that the operator $E^*E$ acts on a 1-form $\theta \in \Omega^1(M)$ according to

$$(E^*E(\theta))_k = (E^*(E\theta))_k = 4(n+1)^2\,(\bar{\Delta}(\bar{\Delta}\theta))_k + \cdots,$$

where the ellipsis denotes lower-order terms. Consequently, the principal symbol of $E^*E$ is given by

$$\sigma_{E^*E}(\xi) = \sigma_{E^*}(\xi)\,\sigma_E(\xi) = 4(n+1)^2\,\|\xi\|^4\,\mathrm{Id}_{T_x^*M},$$

and hence $E^*E$ is elliptic, since its principal symbol is positive definite for $\xi \neq 0$.

We shall call the fourth-order elliptic operator

$$E^*E\colon \Omega^1(M) \to \Omega^1(M)$$

the *Eisenhart operator*.

By virtue of [7, p. 632], for an elliptic operator—and, in particular, for the Eisenhart operator $E^*E$—there holds the $L^2$-orthogonal decomposition of the space of differential 1-forms:

$$\Omega^1(M) = \ker(E^*E) \oplus_{L^2} \mathrm{im}(E^*E).$$

Using standard properties of elliptic operators on a compact (boundaryless) Riemannian manifold $(M, g)$ (see also [7, pp. 631–632]), we have

$$\langle E^*E(\theta), \theta\rangle_{L^2} = \|E\theta\|_{L^2}^2 \geq 0.$$

It follows that the kernel $\ker(E^*E)$ coincides with the kernel of $E$. Therefore, $\ker(E^*E)$ consists of all infinitesimal projective transformations of $(M, g)$, just as does $\ker E$.

Moreover, $\ker(E^*E)$ is finite-dimensional; consequently, $\ker E$ is also finite-dimensional. This agrees with the result cited above on the finite dimensionality of the Lie algebra of infinitesimal projective transformations, which in our setting is precisely the kernel of $E$.

As a result, the corollary to Theorem 1 follows.

**Corollary 5.** *If $(M, g)$ admits no nontrivial infinitesimal projective transformations, $E^*E$ is strictly positive, with*

$$\Omega^1(M) = ker\ E^*E \oplus_{L^2} im\ E^*E.$$

Before stating the next result, recall that a connected compact (boundaryless) $n$-dimensional Riemannian manifold $(M, g)$, $n \geq 2$, of nonpositive sectional curvature is said to have *quasi-negative sectional curvature* if there exists at least one point at which all sectional curvatures are strictly negative (see, e.g., [34], [38]).

**Corollary 6.** *On a connected compact (boundaryless) Riemannian manifold $(M, g)$ of quasi-negative sectional curvature and dimension $n \geq 2$, the Eisenhart operator*

$$E^*E \colon \Omega^1(M) \to \Omega^1(M)$$

*is positive definite.*

**Proof.** Recall that on a compact (boundaryless) manifold $(M, g)$ of dimension $n \geq 2$, the kernel of the second-order linear differential operator

$$E \colon \Omega^1(M) \to C^\infty(T^*M \otimes S_0^2 M),$$

acting according to formula (5.1), coincides with the kernel of the elliptic Eisenhart operator $E^*E \colon \Omega^1(M) \to \Omega^1(M)$.

As in the proof of Theorem 2, define on $(M, g)$ the symmetric tensor field

$$\varphi := 2(n+1)\,\delta^*\theta + 2g(\delta\theta),$$

so that trace$_g\,\varphi = -\delta\theta$. Then equations (5.1) take the form

$$E(\varphi)_{kij} = \nabla_k \varphi_{ij} - \frac{1}{2}(g_{kj}\nabla_i \varphi^l_{\ l} + g_{ki}\nabla_j \varphi^l_{\ l})$$

in local coordinates $x^1, \ldots, x^n$ of an arbitrary chart $(U, \psi)$ on $M$.

Hence, if $\theta \in \ker(E^*E)$, the local components $\varphi_{ij}$ satisfy

$$\nabla_k \varphi_{ij} = \frac{1}{2}(g_{kj}\nabla_i \varphi^l_{\ l} + g_{ki}\nabla_j \varphi^l_{\ l}), \qquad (3.6)$$

from which it follows that

$$\nabla^k \varphi_{kj} = \frac{1}{2}(n+1)\nabla_j \varphi^l_{\ l}.$$

On the other hand, for any connected compact (boundaryless) Riemannian manifold $(M, g)$, the integral identity (3.6) holds for the local components of any symmetric tensor $\varphi \in C^\infty(S^2 M)$.

For a symmetric tensor field $\varphi$ whose components satisfy (3.6), that identity reduces to

$$\int_M \left( \sum_{i<j} \sec(e_i, e_j)(\rho_i - \rho_j)^2 - \frac{1}{4}(n-1)(n+2) \parallel \delta\varphi \parallel^2 \right) dv_g = 0,$$

where $\parallel \delta\varphi \parallel^2 = g(\delta\varphi, \delta\varphi) = \nabla_i \varphi^{ij} \nabla^k \varphi_{kj}$.

It follows that on a compact Riemannian manifold $(M, g)$ of quasi-negative sectional curvature, the only tensor field

$$\varphi \in C^\infty(S^2 M) \cap \ker E$$

is of the form $\varphi = Cg$, where $C$ is constant.

In the present setting this yields

$$\mathscr{L}_\xi g = (2(n+1)^{-1} \text{div } \xi + C) g,$$

whence $\text{div } \xi = \frac{1}{2}nC = \text{const}$. On a compact (boundaryless) Riemannian manifold this implies $\text{div } \xi = C = 0$ (see [35, p. 30]). Therefore, $\mathscr{L}_\xi g = 0$, i.e., $\xi$ is an infinitesimal isometry.

However, it is known (see [35, p. 36]) that no nontrivial isometries exist on a compact (boundaryless) Riemannian manifold satisfying the above curvature condition. Consequently,

$$\ker(E^*E) = \{0\},$$

and the Eisenhart operator $E^*E$ is positive definite.

**Remark.** The result of [33] concerning the dimension of the Lie algebra of infinitesimal projective transformations is local in nature. This is consistent with Corollary 2 proved above, since it excludes the compact hyperbolic space $\mathbb{H}^n$ from consideration.

### §6. Spectral properties of the Eisenhart operator

When studying the spectral properties of the elliptic differential **Eisenhart operator**

$$E^*E: \Omega^1(M) \to \Omega^1(M)$$

on a compact (boundaryless) Riemannian manifold $(M, g)$, we rely on standard material from the monographs [3], [36], and [37] without detailed citation. Using both our previous results and the general theory of elliptic operators, we summarize the spectral properties of $E^*E$, viewed as a self-adjoint, nonnegative elliptic operator.

1. *General spectral properties*

Firstly, the spectrum of the Eisenhart operator $E^*E$ is discrete. Its eigenvalues have finite multiplicity and form an increasing sequence tending to $+\infty$.

If there are no infinitesimal projective transformations on the compact (boundaryless) manifold $(M, g)$, i.e. if $\ker E = \{0\}$, then $E^*E$ has a strictly positive spectral lower bound, since all its eigenvalues are positive:

$$\operatorname{Spec}(E^*E) = \{0 < \mu_1 \leq \mu_2 \leq \cdots \to \infty\}.$$

If ker $E \neq \{0\}$, which may occur in the presence of infinitesimal projective transformations, then 0 is the smallest eigenvalue of $E^*E$ with multiplicity

$$\dim_{\mathbb{R}} \ker (E^*E) \leq n(n+2),$$

and all other eigenvalues are positive:

$$\text{Spec }(E^*E) = \{0 = \mu_0 \leq \mu_1 \leq \cdots \to \infty\}.$$

The eigen–1-forms of the Eisenhart operator form an orthonormal basis in $L^2(\Omega^1(M))$, and the spectral decomposition is compatible with the Hodge decomposition of $\Omega^1(M)$.

2. *The case of the sphere* $(\mathbb{S}^n, g_{\text{can}})$.

We now consider the Eisenhart operator acting on 1-forms on the unit sphere $(\mathbb{S}^n, g_{\text{can}})$, where $g_{\text{can}}$ denotes the canonical metric of constant sectional curvature 1.

Since $(\mathbb{S}^n, g_{\text{can}})$ is simply connected, compact, and of constant curvature 1, the space of harmonic 1-forms is trivial:

$$\ker \Delta_H = \{0\}.$$

Hence the Hodge decomposition reduces to

$$\Omega^1(S^n) = \text{im } d \oplus_{L^2} \text{im } \delta.$$

Every 1-form can be written as

$$\omega = df + \delta\theta,$$

where $df$ is exact and $\delta\theta$ is coexact ($\delta(\delta\theta) = 0$).

The Hodge–de Rham Laplacian $\Delta_H$ acting on 1-forms has spectrum splitting into exact and coexact parts:

- For exact forms $df$:

$$\Delta_H \theta = \lambda_{\text{exact}}(k)\theta, \lambda_{\text{exact}}(k) = k(k+n-1), k = 1,2,\ldots$$

- For coexact forms ($\delta\theta = 0$):

$$\Delta_H \theta = \lambda_{\text{coexact}}(k)\theta, \lambda_{\text{coexact}}(k) = (k+1)(k+n-2), k = 1,2,\ldots$$

3. *Coexact 1-forms*

Let $\theta$ be coexact ($\delta\theta = 0$). Then

$$E\theta = (n+1)\nabla(\mathcal{L}_\xi g),$$

here $\xi$ is the vector field dual to $\theta$.

Suppose

$$E^*E(\theta) = \mu_{\text{coexact}}(k)\theta.$$

Then

$$\langle E^*E\theta, \theta\rangle_{L^2} = \mu_{\text{coexact}}(k) \parallel \theta \parallel_{L^2}^2.$$

Since $E^*E$ is self-adjoint,

$$\langle E^*E\theta, \theta\rangle_{L^2} = \langle E\theta, E\theta\rangle_{L^2} = (n+1)^2 \langle \nabla(\mathcal{L}_\xi g), \nabla(\mathcal{L}_\xi g)\rangle_{L^2}.$$

Using curvature identities on the sphere and the Weitzenböck formula

$$\Delta_H \theta = \bar{\Delta}\theta + \text{Ric}(\xi, \cdot),$$

with Ric $= (n-1)g$, we obtain

$$\langle E^*E\theta, \theta\rangle_{L^2} = (n+1)^2 (\langle \Delta_H \theta, \theta\rangle_{L^2} + \parallel \theta \parallel_{L^2}^2).$$

Hence,

$$\mu_{\text{coexact}}(k) = (n+1)^2 \left( \frac{\langle \Delta_H \theta, \theta\rangle}{\parallel \theta \parallel^2} + 1 \right).$$

If $\theta$ is an eigenform of $\Delta_H$,

$$\Delta_H \theta = \lambda_{\text{coexact}}(k)\theta,$$

Then

$$\mu_{\text{coexact}}(k) = (n+1)^2((k+1)(k+n-2)+1).$$

## 4. Exact 1-forms

Let $\theta = df$ and suppose
$$E^*E(df) = \mu_{\text{exact}}(k)\, df.$$
Using the well-known formula
$$(E\,\theta)_{kij} = 2(n+1)\,\nabla_k(\delta^*\theta)_{ij} + 2\,g_{ij}\,\nabla_k(\delta\theta) + + g_{kj}\,\nabla_i(\delta\theta) + g_{ki}\,\nabla_j(\delta\theta),$$
we obtain
$$(df)_{kij} = 2(n+1)\nabla_k(\text{Hess }f)_{ij} + 2g_{ij}\nabla_k(\bar{\Delta}f) + g_{kj}\nabla_i(\bar{\Delta}f) + g_{ki}\nabla_j(\bar{\Delta}f), \qquad ($$
where Hess $f = \nabla df$.

Employing integration by parts and the Ricci identity on the sphere,
$$\bar{\Delta}(\nabla f) = \nabla(\bar{\Delta}f) + (n-1)\nabla f,$$
we obtain
$$E^*E(df) = (n+1)^2(d(\bar{\Delta}f) + 2n\, df) = (n+1)^2(\bar{\Delta} + 2n)\, df.$$
If $f$ is a spherical harmonic of degree $k$,
$$\bar{\Delta}f = \lambda_{\text{scalar}}(k)f, \lambda_{\text{scalar}}(k) = k(k+n-1),$$
then for $k \geq 2$,
$$E^*E(df) = (n+1)^2(k(k+n-1) + 2n)\, df.$$
Thus,
$$\mu_{\text{exact}}(k) = (n+1)^2(k(k+n-1) + 2n).$$
Therefore, on the sphere $(\mathbb{S}^n, g_{\text{can}})$, the spectrum of the Eisenhart operator splits into exact and coexact parts and can be computed explicitly in terms of classical spherical harmonic eigenvalues.

## Conclusion

In this paper we introduced and systematically investigated two natural two-order and fourth-order elliptic operators arising in projective Riemannian geometry: the *Sinjukov operator* $S^*S$, acting on symmetric 2-tensors, and the *Eisenhart operator* $E^*E$, acting on differential 1-forms, respectively. Both operators were constructed as compositions of first-order overdetermined operators with their formal adjoints and were studied within the framework of global elliptic theory on compact Riemannian manifolds.

We established that the operators $S^*S$ and $E^*E$ are elliptic, self-adjoint, and nonnegative on compact manifolds without boundary. Their principal symbols were shown to be injective, which implies overdetermined ellipticity of the underlying first-order operators and guarantees finite-dimensional kernels. As a consequence, we obtained $L^2$-orthogonal decompositions of the corresponding function spaces and clarified the geometric meaning of the kernels:

- $\ker(S^*S) = \ker S$ consists of symmetric 2-tensors generating projectively equivalent metrics;
- $\ker(E^*E) = \ker E$ coincides with the Lie algebra of infinitesimal projective transformations.

Under curvature assumptions, these kernels were described explicitly. In particular, on compact manifolds of nonpositive sectional curvature, elements of $\ker(S^*S)$ are parallel symmetric tensors, while under quasi-negative curvature the kernel reduces to multiples of the metric. Analogously, for quasi-negative sectional curvature the Eisenhart operator is strictly positive, implying the absence of nontrivial infinitesimal projective transformations.

We also analyzed the spectral structure of both operators. Using general elliptic theory, we proved discreteness of the spectrum, finite multiplicity of eigenvalues, and compatibility with Hodge-type decompositions. In the model case of the standard sphere $(\mathbb{S}^n, g_{\text{can}})$, the spectra of $S^*S$ and $E^*E$ were computed explicitly by decomposing tensor and form spaces into irreducible components and expressing eigenvalues in terms of the scalar, vector, and Lichnerowicz Laplacians. These computations demonstrate how the geometry of constant curvature determines the full spectral behavior of the operators.

Overall, the results clarify the analytic structure underlying projective transformations of Riemannian manifolds and provide a spectral framework for further study of rigidity, curvature restrictions, and geometric invariants associated with projective equivalence.

**Competing interests**

The authors declare that they have no competing interests.

**Author's contributions**

All authors contributed equally to the writing of this paper. All authors read and approved the final manuscript.

**Josef Mikeš**

Dept. of Algebra and Geometry, Palacky University,
771, 80 Olomouc, Czech Republic
Mikes@risc.upol.cz

**Sergey E. Stepanov**

Department of Mathematics and Data Analysis, Finance University under the Government of Russian Federation,  49-55, Leningradsky Prospect, 125468 Moscow, Russia;
Department of Mathematics, Russian Institute for Scientific and Technical Information of the Russian Academy of Sciences, 20, Usievicha street, 125190 Moscow, Russia;
s.e.stepanov@mail.ru